# A topological approach to creating any *pulli kolam*, an artform from South India

Venkatraman GOPALAN*, Brian K. VANLEEUWEN

[1] *Materials Science and Engineering, Pennsylvania State University, University Park, PA 16802, USA*

*E-mail address: vgopalan@psu.edu*

**Abstract:** *Pulli kolam* is a ubiquitous art form in south India. It involves drawing a line looped around a collection of dots (*pulli*s) place on a plane such that three mandatory rules are followed: all line orbits should be closed, all dots are encircled and no two lines can overlap over a finite length. The mathematical foundation for this art form has attracted attention over the years. In this work, we propose a simple 5-step topological method by which one can systematically draw all possible *kolam*s for any number of dots *N* arranged in any spatial configuration on a surface.

## 1. What is a *Kolam*?

Figure 1 depicts an example of a *kolam*, an ancient and still popular South Indian art form. This particular type of *kolam* is called the *pulli kolam* in Tamil, which consists of a series of dots (called *pullis*) placed on a surface, each of which is then circumscribed by lines that form closed orbits. It is a very common sight on the threshold of homes in the five southern states with a combined current population of ~252 million. They are called by varied names in the respective regional languages of these states: *kolam* in Tamil spoken in Tamil Nadu, *golam* in Malayalam spoken in Kerela, *rangole* in Kannada spoken in Karnataka, and *muggulu* in Telugu spoken in Andhra Pradesh and Telangana. With every sunrise, women wash the floor in front of the houses, and using rice flour, place the dots and draw a *kolam* largely from memory. Learning how to draw *kolam*s from an early age is an important aspect of growing up in southern India, especially for girls. As they continue to learn from other women in their family, the *kolam*s become increasingly complex, with a larger number of dots and more intricate line orbits. Remembering the dot configurations and line orbits is a daily exercise in geometric thinking. The process is immensely pleasurable, especially when a *kolam* is successfully completed with no loose ends.

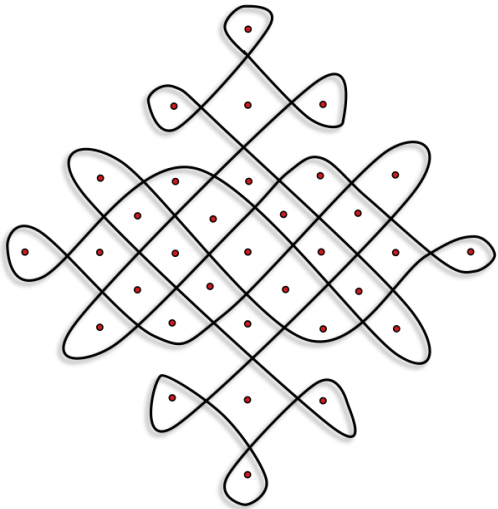

Fig. 1. Example of a *pulli kolam* called Brahma's knot.

While the conventional *kolams* impose several constraints, here we begin with three simple rules in order to give ourselves greater room for discovery and creativity. Given an arbitrary arrangement of dots on a plane, the following three *mandatory* (**M**) *constraints* define a *kolam*:

**M**1: All dots should be circumscribed.

**M**2: All interactions between two lines must be at points, i.e. two line segments cannot overlap over a finite length.

**M**3: All line orbits should be closed, i.e. no loose ends.

In addition to the above rules, one may choose to apply additional *optional guidelines* (**O**). There is no limit to the number of such optional guidelines that can be followed, but we will explore some of them later in this work.

While *kolam*s are widely rendered from memory, the process of creating entirely new ones, especially complex *kolam*s with a large number of dots is far more challenging. This work attempts to provide a simple 5-step method by which anyone can create a very large number of *kolam*s from any arbitrary pattern of dots. The proposed topological method deemphasizes memory; in principle, anyone who knows just the method will be able to draw a large number of

*kolam*s with no other prior knowledge.

Many previous pioneering works exist that have provided mathematical insights into the form of a *kolam* over the past four decades. These include converting *kolam*s into numbers and linear diagrams [1], using graph, picture, and array grammers [2-10], extended pasting schemes [11], morphism of monoids [12], L- and P-systems [13,14], gestural lexicons [15], knot theory [16], and mirror curves [17]. Of these, the work of Yanagisawa & Nagata [1], has similarities to this work. They begin with 5 rules for *kolam*, define square unit tiles that can be assembled into larger *kolam*s, define two types of nearest neighbor interactions between dots (line crossing, 1, or uncrossing, 0) and convert these tiles into binary number arrays. Nagata [18] also addressed the construction of a primitive *kolam* for an arbitrary dot array with a similar approach. In contrast, the work presented here has a purely topological approach: it defines only 3 mandatory rules for defining a *kolam*, has no standard tiles, generalizes the ideas to any arbitrary arrangement of dots arranged in any shape (not necessarily square arrays), generalizes to interactions between any two dots (instead of only the nearest or next nearest neighbors), and to three or more number of bonds between an interacting pair of dots. The work suggests that for a given number of dots, $N$, there are a limited number of *parent kolam types* from which all other *kolam*s originate. All parent *kolam*s within a parent *kolam type* are homotopic (or topologically equivalent).

## 2. How many *kolam*s for one dot ($N$=1)?

Figure 2 depicts a single dot, and a variety of lines circumscribing it that follow the three mandatory rules mentioned above. The *kolam* in general could be amorphous in shape, as in Fig 2a, and in the special case of Fig. 2b is a circle. Multiple circumscriptions around the dot are possible, as in Figs 2c, d, and e.

It becomes immediately clear from Fig. 2 that the number of possible *kolam*s thus defined, with only the mandatory rules, is *infinite*. One may arbitrarily impose additional optional guidelines (**O**) to limit the number of *kolam*s. Here are some:

**O**1: Only one circumscription of the line is allowed around each dot.

**O**2: A line circumscribing a dot should be as resourceful (simple) as possible, without additional unnecessary wiggles or flourishes (e.g. Fig. 2b is resourceful vs. Fig. 2a is not).

**O**3: While a *kolam* may be created by a minimum number of dots $N$ needed for the 5-step method proposed below, one can then eliminate dots from, or add dots to, or move dots in a *kolam after* it has been drawn, provided the process does not violate the mandatory rules. The final *kolam* may thus appear to have $N_{final}$ dots, where $N_{final}$ may or may not be equal to $N$.

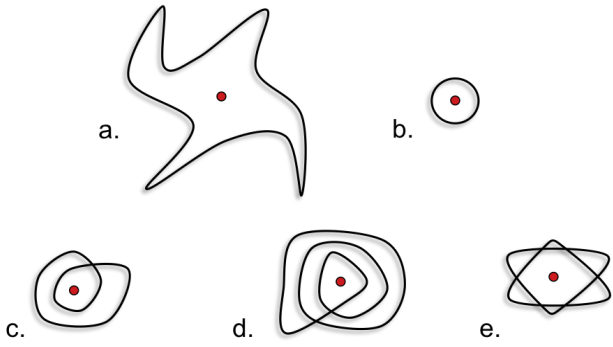

Fig. 2. Examples of *kolam*s around one dot that follow the mandatory rules **M**1- **M**3. An infinite number of *kolam*s are possible. Additional optional (**O**) guidelines, **O**1-**O**6, can limit the number allowed.

With **O**1 restriction, only 2a and 2b survive. With **O**1 and **O**2, only 2b will survive. Figure 2e, depicting a Star of David is a common *kolam*, which apparently is eliminated by **O**1. However, this *kolam* can also be generated by placing six dots ($N=6$), one inside each ray of the star, and following the 5-step method proposed below. The 6 dots may later be erased, and one dot placed in the middle ($N_{final}=1$) according to **O**3 to generate Fig. 2e. Another example is the Brahma's knot in Fig. 1, which can be generated by only $N=25$ dots. However, Fig. 1 has $N_{final}=33$ dots; the additional two horizontal rows of 4 dots each (total of 8 dots) in that *kolam* would be placed (according to **O**3) after constructing the *kolam* with only 25 dots by the method proposed below.

## 3. Method to construct *kolam*s for an arbitrary arrangement of $N_{min}$ dots

First, we define several types of bonds (*b*) between a pair of dots, as shown in Fig. 3. The X- and the B- bonds were discussed in Yanagisawa & Nagata [1] and they were indexed as a line crossing, 1, or an uncrossing, 0. The D-bond corresponds to additional variation (a type of two-dot joining, indexed as 2) over the pictorial code proposed by Nagata [18]. In general, there are infinitely many possible bond types, such as 2X-, $X^+$- and $X^-$- bonds and so on. From here on, we will focus on the cross (X)-bond, the double (D)-bond, and the broken (B)-bonds ($b=3$) in this work.

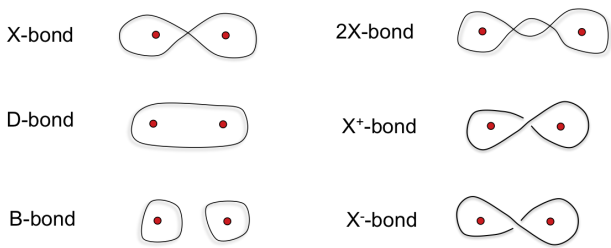

Fig. 3. Infinitely many types of bonds are possible between a pair of dots that follow the rules **M**1, **M**2, **M**3, and guidelines **O**1 and **O**2, some of which are shown here. The $X^+$ and $X^-$ bonds have one arm of the cross on top of the other arm, thus utilizing the third dimension out of the page. They are shown here as a broken bottom arm for ease of depiction on a plane.

Next, we propose a 5-step method to build all possible *kolam*s for an *arbitrary* pattern of $N$ dots in two-dimensions. These rules are illustrated for a simple 2-dot ($N=2$) case in Figure 4.

**Step 1**: Place the dots in *any* configuration of your choice in 2-dimensions.

**Step 2**: Draw a perpendicular bisector line segment between every pair of dots in the general case of following only rules **M**1-**M**3. (More generally, this line segment does not need to be a bisector, and does not need to lie between the two dots.) This bisector has analogy to the navigating line (N-line) used by Nagata [18].

**Step 3**: Draw closed ghost-like figures around each dot, which we will playfully call *squishies*, suggesting that they are freely deformable. There will be $N$ squishies for $N$ dots. Each squishy will have $J$ arms that touch a corresponding arm from a different squishy pairwise at the bisector line, leading to $J$ junctions. We will call this structure, the *parent kolam*. All *kolam*s arising from the $N$ dots will arise from this parent *kolam*.

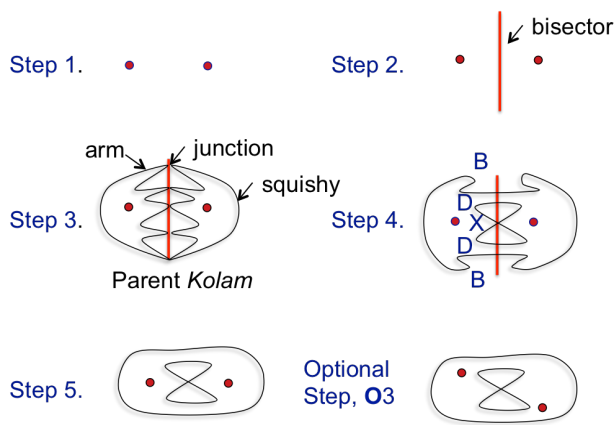

Fig. 4: Illustrating the construction of a *kolam* in 5 steps plus optional guideline **O**3: The procedure is shown for $N=2$ (2 dots) and $J=5$ (5 junctions). If each junction is restricted have one of 3 types of bonds (X-, D-, or B-), it can lead to $3^5$ =243 possible *kolam*s. One of these *kolam*s, namely, B-D-X-D-B, is shown in the figure in Step 4. In the optional guideline, the dots have been rearranged as an example of rule **O**3 after the *kolam* is drawn in Step 5.

**Step 4**: Now start drawing the *kolam* from any point on a squishy, and follow along until you reach a junction. Then transform that junction into a *cross-bond* (X-bond), a *double-bond* (D-bond), or a *broken-bond* (B-bond). Continue in a similar way until you return to the starting point. If some dots are still not encircled, start a new line from a squishy around one of those remaining dots, and continue till you return back to the start of that line. Repeat this process till the *kolam* is complete and all the dots are encircled.

**Step 5:** Smooth the curves so that the lines are resourceful according to **O**2. This will result in a *kolam* that *will* obey the rules **M**1, **M**2, and **M**3.

As an optional guideline, you can eliminate any or all dots, or add new dots, or move the existing dots according to **O**3. In addition, one may impose further optional guidelines to whittle down the number of *kolam*s:

**O**4: Only the nearest neighbor dots interact through bonds other than broken bonds. All other bonds are broken.

**O5**: Only one junction ($J=1$) is allowed between one pair of dots.

**O6**: Symmetry equivalent junctions in the parent *kolam* will have the same type of bonds. To find sets of symmetry equivalent junctions, visual inspection of possible rotations axes and mirror symmetries is recommended. For a mathematical approach, find the point group of the arrangement of dots, and using the symmetry operations of the point group, see which set of junctions transform into each other.

In general, with rules **M**1, **M**2, **M**3 and optional guidelines **O**1 & **O**2 in place, with $J$ number of junctions per pair of dots, $N$ and with $b$ types of bonds allowed (Fig. 3), one can write the number, $K$, of possible *kolams* as:

$$\# Kolams = K = b^{JN(N-1)/2} \qquad (1)$$

where the exponent of $b$ is the number of possible junctions between all possible pairs of dots. For example, if $N=2$ (2 dots), $J=1$ (1 junction) and $b=3$ (3 bonds), then $K=3$. These 3 *kolams* are shown in Figure 3. Obviously, $K$ gets large very quickly as $J$, $b$ and $N$ increase. In the rest of this work, we will restrict ourselves to $J=1$ and $b=3$.

If the optional guideline **O**6 is imposed in addition, and symmetry equivalent junctions identified, let there be $g$ groups, each containing $S_g$ number of symmetry equivalent junctions, such that $\sum_g S_g = JN(N-1)/2$. Then the number of possible *Kolams* (Eq. 1) can be revised as $K = b^g$.

Note that we assert in step 5 that this procedure will always result in a *kolam* that obeys the mandatory rules. This arises from the rules of construction. The parent *kolam* is always drawn in the above steps in such a way as to not violate the three mandatory (**M**) rules: all dots are circumscribed by squishies and there are no loose ends in the parent *kolam*. Nor does the transformation of the junctions in Step 3 violate this rule: the bonds where lines cross, e.g. the X- bond, cross at a single point per crossing. Hence the final *kolam* also follows the minimal mandatory rules **M**1-**M**3.

## 4. Exploring *Kolams* with 3 dots ($N=3$).

The number of possible *kolams* for $N=3$ following rules **M**1-**M**3 and optional guidelines **O**1, **O**2, and **O**5 ($J=1$) can be computed from Eq. 1 as $K=3^{1 \times 3 \times (3-1)/2} = 3^3 = 27$. Two different *parent kolams* for $N=3$ are shown in Fig. 5.

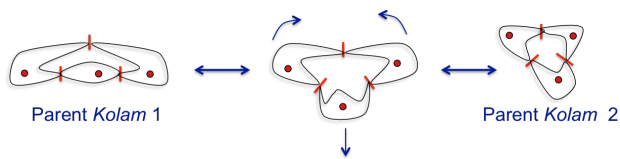

Fig.5. Two possible parent *kolams* for 3-dots ($N=3$) and $J=1$. The intermediate structure shows how one can distort parent 1 into parent 2, demonstrating that they are homotopic.

Parent 1 places the three dots on a line, while parent 2 places them in a triangle. These two parent *kolams* are topologically equivalent, or *homotopic*. In other words, a continuous distortion of one structure can result in the other without cutting or breaking bonds, as shown by a transformation through the intermediate structure in Fig. 5. Hence, every one of the 27 *kolam*s derived from parent *kolam* 1 will have a topologically equivalent cousin *kolam* derived from parent *kolam* 2. Thus we can conclude that for $J=1$, all $N=3$ *kolam*s arise from a single parent *kolam type*.

The 27 *kolam*s derived from parent *kolam* 2, with the special case of the 3 dots arranged in an equilateral triangle, are shown in Fig. 6.

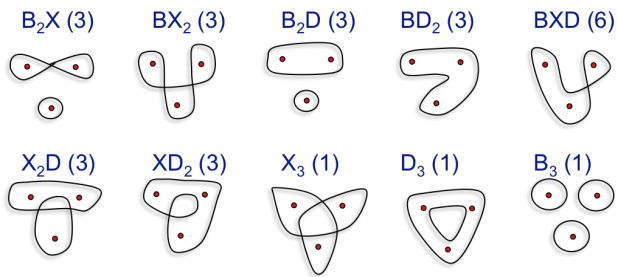

Fig. 6: The 27 *kolam*s generated from 3 dots ($N=3$) and $J=1$. There are 3 possible pairs of dots. The notation, $B_2X$ (3), for example indicates that two of the pairs have broken-bonds and one pair has a cross-bond. The (3) in the end indicates that three such *kolam*s of the same *type* exist, generated by the permutation of the X-bond between the three pairs in the case of $B_2X$.

Did we find all possible *kolam*s with $N=3$? If so, how about the *kolam* on the left in Fig. 7? It turns out that this *kolam* is captured by the proposed method for $N=4$, where an additional dot is placed in the middle of Figure 7. This is discussed in the next section. The example is again illustrative of the fact that a *kolam*, once created, is distinctive in its own right, irrespective of the presence or absence of dots. *The characteristic N for a given kolam may be defined as the minimum number of dots required for generating the kolam with the above 5-step method.* However, note that when dots are removed or added to a *kolam*, the resultant *kolam*s may no longer be topologically equivalent to the original *kolam*.

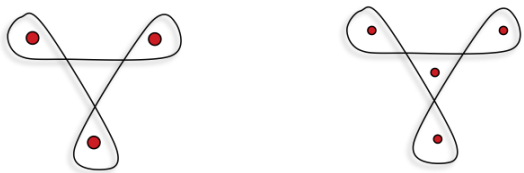

Fig. 7: The above *kolam* on the left might appear to be a 3 dot *kolam*. However, it is an $N=4$ dot *kolam* (above right) created using the 5-step method described above. By erasing the center dot in the right *kolam*, one can generate the $N_{final}=3$ *kolam* according to **O**3. These two *kolam*s are *not* homotopic.

## 5. Exploring *Kolam*s with 4 dots ($N=4$).

Three different configurations of parent *kolam*s are shown in Fig. 8 for $N=4$.

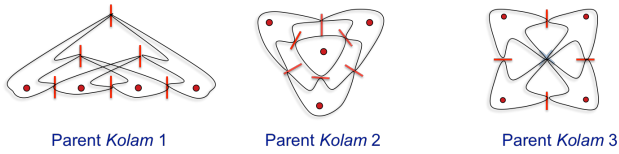

Parent *Kolam* 1   Parent *Kolam* 2   Parent *Kolam* 3

Fig. 8: Three types of parent *kolam*s for $N=4$ dots and $J=1$. Parents 1 and 3 are homotopic and form one parent *type*. Parent 2 forms a second parent *type*.

It is possible to show that parents 1 and 3 are homotopic. Such equivalence is shown in Fig. 9a, and hence they form a single parent *type*. However, parent 2 forms a distinct parent type as shown in Fig. 9b. The number of possible *kolam*s for any parent *kolam* with $N=4$ following rules **M1**-**M3** and **O1**, **O2**, and **O5** can be computed from Eq. 1 as $K=3^{1\times 4 \times (4-1)/2} = 3^6 = 729$.

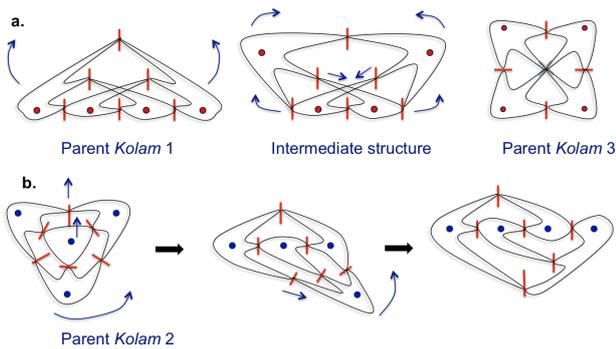

a. Parent *Kolam* 1   Intermediate structure   Parent *Kolam* 3
b. Parent *Kolam* 2

Fig. 9. (a) Parent *kolams* 1 and 3 for $N=4$ and $J=1$ shown in Figure 8 are demonstrated to be topologically equivalent by continuously deforming parent 1 into 3 in panel a; hence they form a single parent *type*. Panel b shows that distorting parent 2 in Figure 8 does not lead to parent 1; hence they are distinct parent types.

The 729 possible *kolam*s from each parent is a large number, and so we choose here to impose additional restrictions in order to explore only a subset. For example, optional guideline **O6** suggests that symmetry equivalent junctions will have the same type of bond.

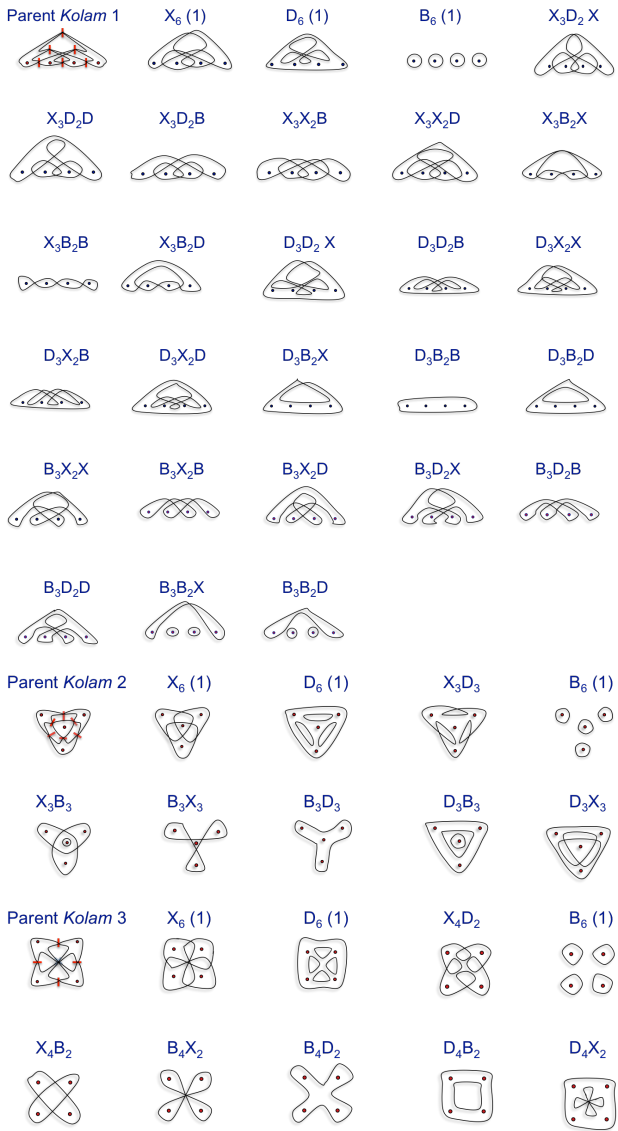

Fig. 10: Four dot (*N*=4) *kolam*s derived from the three parents in Fig. 8, under the constraints of the rules **M**1-**M**3, and optional guidelines **O**1(only one circumscription per dot), **O**2 (simplifying the line), **O**5 (*J*=1), and **O**6 (symmetry equivalent junctions will have the same type of bond). Note that the parent kolams in this figure have been chosen in the special shapes of a line (parent 1), an equilateral triangle (parent 2) and a square (parent 3). These choices as well as the optional guidelines eliminate many kolams that a reader might otherwise be able to visualize.

This allows for the symmetry of the parent phase to be preserved while bonds are formed. The various *kolam*s derived from three different parent *kolam*s (1, 2, and 3) in Fig. 8 under the rules of **M**1-**M**3 and **O**1, **O**2, **O**5, and **O**6 are shown in Fig. 10. For parent *Kolam* 1 in Figure 10, there are 3 groups (*g* =3) of symmetry equivalent junctions related by a vertical mirror symmetry. Thus the number of *Kolams* with *J*=1 is $K=3^3=27$. For both the special cases of parent *Kolam* 2 (dots forming an equilateral triangle) and *Kolam* 3 (dots forming a square), *g*=2 arising from a 3-fold and 4-fold rotational axes respectively, and hence $K=3^2=9$ as shown. We note that $B_3X_3$ with *N*=4 captures the *kolam* that was missed in Fig. 7 by *N*=3.

## 6. Conclusions

We have demonstrated a method of generating countless *kolam*s from user-defined dot arrangement on a surface. This method can be mastered by anyone without the need to understand the detailed mathematics behind *kolam*s. For a give number, *N*, the number of possible *kolam*s that follow only the mandatory rules **M**1-3 is infinite, even for a 1-dot *kolam* (*N*=1). However, by following additional optional guidelines **O**1 and **O**2, this number is finite as given by Eq. 1. Addition of guideline **O**6 modifies this equation.

We show by example that for a given number of dots *N*, a set of parent *kolam types* exist, from which all possible *kolam*s can be generated. All parent *kolam*s within a single *type* are homotopic. Hence the resultant *kolam*s from these homotopic parents will also have corresponding homotopic cousins. Though a rigorous proof for such homotopy in general has not been presented, it can be argued based on the method of construction similar to that shown in Fig. 9.

Kolams with higher *N* get richer and more complicated quickly. For example, Figure 11 shows an example parent *kolam* for *N*=1 and *J*=1, and two possible children *kolam* arising from it. The readers are encouraged to try generating other parent and children *kolams* for this case.

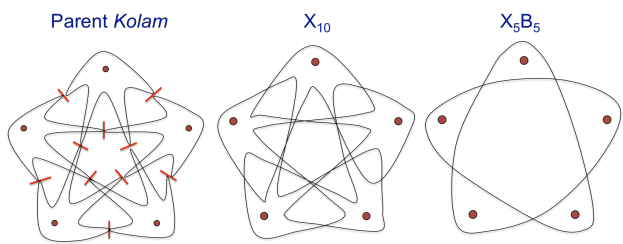

Figure 11: An example of a parent *kolam* and two children *kolams* for *N*=5 and *J*=1.

There are several advantages to this simple method: (1) It is applicable for any number of dots, *N*, arranged in *any* configuration in 2-dimensions. (2) While the proposed method may not always guarantee aesthetics, it is simple enough for a user to impose additional aesthetically appropriate optional guidelines. (3) A computer program can vary *b*, *J*, and *N* for generating numerous *kolam*s following the three mandatory rules, plus any number of user-defined optional guidelines. This leads to the possibility of creating an interactive website or a mobile app that can help a user to generate *kolam*s at will. Such an app will get the user involved in the creative process, including young children who may be introduced to art, symmetry and topology through *kolam*s. The method is also applicable to generating other similar patterns such as some of the Chinese knots by using two-level crossing bonds such as $X^+$ and $X^-$ bonds in Fig. 3.

**7. Acknowledgments:** The authors would like to thank Mikael Rechtsman and Chaoxing Liu from the Pennsylvania State University for exciting and valuable discussions. The ideas for this work tangentially originated from developing the theory of distortions funded by the National Science Foundation grant number DMR-1420620 and DMR-1210588.